\documentclass[epsbox,12pt,reqno]{amsart}

\usepackage{amsmath,amssymb,latexsym,amscd,amsthm} 
\usepackage{graphicx}

\newtheorem{theorem}{Theorem}[section]  
\newtheorem{corollary}[theorem]{Corollary}
\newtheorem{proposition}[theorem]{Proposition}
\newtheorem{lemma}[theorem]{Lemma}
\theoremstyle{definition}

\newtheorem{question}[theorem]{Question}
\newtheorem{conjecture}[theorem]{Conjecture}
\theoremstyle{remark}
\newtheorem{remark}[theorem]{Remark}

\def\sd{\operatorname{depth}}
\def\csd{\mathbf{depth}}
\def\S{\mathbf{S}}

\def\T{\mathbf{S}''}

\def\f{\mathbf{f}}
\def\g{\mathbf{g}}
\def\core{\operatorname{core}}
\def\cone{\operatorname{cone}}
\def\conv{\operatorname{conv}}
\def\aff{\operatorname{Aff}}
\def\south{\operatorname{south}}
\def\north{\operatorname{north}}
\def\interior{\operatorname{int}}
\def\bara{B\'ar\'any}
\def\cara{Carath\'eodory}
\def\R{\mathbb{R}}

\def\X{\text{X}}
\def\P{\text{H}}
\def\Sph{\mathbb{S}}
\def\B{\mathbb{B}}
\def\m{\operatorname{m}}
\def\zero{{\bf 0}}
\def\Cap{\operatorname{Cap}}
\def\Can{\operatorname{Can}}
\def\Pole{\operatorname{Pole}}

\title{Colourful Simplicial Depth}

\author{Antoine Deza, Sui Huang, Tamon Stephen and Tam{\' a}s Terlaky}

\address{Advanced Optimization Laboratory,
Department of Computing and Software,
Faculty of Engineering,
1280 Main St.~West,
McMaster University,
Hamilton, Ontario,
Canada
L8S 4K1.
}

\email{\{deza,huangs3\}@mcmaster.ca, tamon@optlab.mcmaster.ca, terlaky@mcmaster.ca}

\thanks{Research supported
by NSERC Discovery grants for the four authors,
by the Canada Research Chair program for the first and last authors
and by a MITACS grant for the second and third authors.}


\begin{document}

\def\minus{-}
\setcounter{page}{1}

\maketitle

\begin{abstract}
Inspired by \bara's colourful \cara~theorem \cite{Bar82}, we introduce
a colourful generalization of Liu's simplicial depth \cite{Liu90}.
We prove a parity property and conjecture that the minimum
colourful simplicial depth of any core point in any $d$-dimensional
configuration is $d^2+1$ and that the maximum is $d^{d+1}+1$.  
We exhibit configurations attaining each of these depths,
and apply our results to the problem of bounding
monochrome (non-colourful) simplicial depth.
\end{abstract}

%
%
\section{Introduction}\label{se:intro}

In statistics there are several measures of the depth of a
point $p$ in ${\R}^d$ relative to a fixed set $S$ of sample points. 
Two recent surveys on data depth are \cite{Alo05} and \cite{FR05},
see references therein.  
The depth measure we are interested in is the {\it simplicial depth} 
of $p$, which is the number of simplices generated by points
in $S$ that contain $p$.
A point of maximum simplicial depth can be viewed as a type of
$d$-dimensional median.  We would like to obtain a lower bound for 
the depth of simplicial medians.

To do this, we consider a generalized
problem where the sample points are colourful.  That is, in
dimension $d$ we consider sample points given in each of at least
$(d+1)$ colours.  Then we define the {\it colourful simplicial depth} of a
point $p$ relative to this sample to be the number of {\it colourful
simplices} (i.e.~simplices with one vertex of each colour) that
contain $p$.  We focus on the situation where the point $p$ is
in the intersection of the convex hulls of the individual colours,
which is called the {\it core} of the configuration.

If $p$ is a core point we would typically
expect the simplicial depth of $p$ to be more than exponential
in $d$.  However, we exhibit configurations where $p$ is a core point
but is contained in only $d^2+1$ colourful simplices.  We
conjecture that any core point $p$ of any $d$-dimensional 
colourful configuration is contained in at least $d^2+1$ 
colourful simplices.  
Along the way, we notice that both in the colourful and 
monochrome cases the simplicial depth of points in general position
(relative to the sample set) sometimes has pleasant parity properties.
We conclude by mentioning some other natural 
problems relating to the colourful
and monochrome simplicial depth.

%
%
\section{Definitions and Background}\label{se:defs}

\subsection{Simplicial Depth}\label{se:depth}
The (closed) simplicial depth of a point $p$ relative to a set 
$S$ of $n=|S|$ points in ${\R}^d$ is the number of (closed) simplices
generated by sets of $(d+1)$ points from $S$ containing $p$ 
in their convex hull.
This was introduced by Liu \cite{Liu90} as a measure of how 
representative $p$ is of the points in $S$.
Denote the simplicial depth of $p$ relative to $S$ as $\sd_S(p)$.
The simplicial depth of $p$ can be interpreted
as the probability that $p$ lies in a random simplex of $S$ 
times a constant factor of $n^{d+1}$ if we sample points from $S$
uniformly with replacement, or times $\binom{n}{d+1}$ if we 
sample without replacement.

We are most interested in the case when $S \cup \{p\}$ is in 
{\it general position}, that is for all $k < d$ there are 
no $k$-dimensional affine 
subspaces contain $k+2$ points from $S \cup \{p\}$.
With this assumption, $p$ will always be in the
interior of any simplices that contain it, so the notions of closed 
and open simplicial depth coincide.  
Without this assumption the closed simplicial depth will be larger.

For a set of points $S$, define $f(S)$ to be the maximum simplicial
depth of a point $p$ relative to $S$, that is:
\begin{equation}\label{eq:maxf}
f(S) = \max_{p \in R^d} \sd_S(p)
\end{equation}
A point $p$ maximizing $f(S)$ can be understood as a 
higher dimensional median point.  We will call any such $p$
a {\it simplicial median}.  Indeed for $d=1$, this is the
usual definition of a median ${\R}$. In higher dimensions,
this definition retains many
desirable properties of the median, such as affine
invariance and a high breakdown point (see e.g.~\cite{Alo05}, 
\cite{FR05}, \cite{GSW92}, \cite{Liu90}).
However, this maximum will not be attained at a point in general position.

We will consider a similar quantity, the maximum simplicial
depth of a point $p$ that maintains $S \cup \{p\}$ 
in general position:
\begin{equation}\label{eq:maxg}
g(S) = \max_{S \cup \{p\} \text{ in general position}} \sd_S(p)
\end{equation}
Equivalently, $g$ is the maximum open simplicial depth of a point
${\R}^d$.  In this way the definition of $g$ can be extended to 
the case when $S$ is not in general position.
While the maximum in (\ref{eq:maxf}) will be 
attained on a discrete set of points in ${\R}^d$, the
maximum in (\ref{eq:maxg}) will be attained on an open set.
For non-empty $S$, we will have $g(S)<f(S)$.

\subsection{Colourful Simplicial Depth}\label{se:csd}
Now consider a situation where points are given in each of $r \ge d+1$
colours.  Then the sample consists of colourful sets 
$S_1, S_2, \ldots, S_r$ which define a colourful configuration $\S$.  
In the following we use a bold font for colourful objects.
A {\it colourful simplex} from these sets
is any simplex whose vertices are chosen from distinct sets.
We define $\csd_\S(p)$ the {\it colourful simplicial depth} of $p$ 
relative to the configuration $\S$ as the number of colourful simplices 
containing $p$.  As with monochrome simplicial depth, colourful
depth can be interpreted probabilistically.  In the case where 
$r = d+1$, colourful depth corresponds to specifying separate 
distributions for each vertex of the simplex.
Dividing the depth by $|S_1| \cdot |S_2| \cdot \ldots  \cdot|S_{d+1}|$ 
gives the probability that $p$ lies in a random colourful simplex
(sampled uniformly).

A choice of sets $S_1, \ldots, S_r$ specifies a {\it colourful
configuration} $\S$ of points.  We call the intersection of the 
convex hulls of the $S_i$'s in a configuration
the {\it core} of $\S$.  \bara~proved
that core points are contained in some colourful simplex; 
this is known at the Colourful \cara~ Theorem \cite{Bar82}.
In the remainder of the paper, 
except where noted, we assume that all configurations and $p$ are in general
position and have a non-empty (hence full-dimensional) core.  
We remark that our results hold under weaker conditions, such as
$p$ not lying on any hyperplanes generated by points from the
configuration.

\subsection{Background}\label{se:background}
Even before the notion of simplicial depth was studied in statistics,
the question of computing bounds for $f(S)$ and $g(S)$ given $n$ and $d$
was studied in
the combinatorics and computational geometry communities.
The two-dimensional question dates back at least to K{\'a}rteszi \cite{Kar55}
who showed that for $n$ points in the plane, $g(S)$ is at most
$\frac{1}{24}(n^3-n)$ for odd $n$ and at most $\frac{1}{24}(n^3-4n)$ 
for even $n$, 
and showed that these bounds were attained when $S$ is the set of
vertices of a regular $n$-gon.
In the early 1980's, Boros and F{\"u}redi \cite{BF84} 
showed $g(S)$ is at least $\frac{n^3}{27} + O(n^2)$, and gave configurations
achieving this bound.

In a beautiful paper, \bara~\cite{Bar82} gave bounds for the 
monochrome simplicial depth in dimension $d$ as an application 
of his Colourful \cara~Theorem.  He obtained a lower bound by
showing that after colouring the points, some point $p$ must be
contained in many colourful simplices.  
A key point of \bara's proof is that a core point $p$ of
a colourful configuration must lie in at least {\it one} colourful
simplex.  Using this fact, for a set $S$ of $n$ points in general
position in ${\R}^d$ \bara~obtains a lower bound of:
\begin{equation}\label{eq:lb}
g(S) \ge \frac{1}{(d+1)^{d+1}}\binom{n}{d+1}+O(n^d)
\end{equation}
This result is asymptotically sharp up to a constant factor as function
of $n$ (for fixed $d$).  However, as \bara~remarks, the constant is
probably quite far from the truth.  Indeed, he gives a sharp upper
bound of:
\begin{equation}\label{eq:ub}
g(S) \le \frac{1}{2^d (d+1)!} n^{d+1} + O(n^d)
\end{equation}
We speculate that the true lower bound is not much less than the 
upper bound.  

One way to improve (\ref{eq:lb}) would be to show that a core point $p$
must lie in more than one colorful simplex.  In \bara's original
paper, he notes that $p$ must in fact lie in at least $(d+1)$ colourful 
simplices, thereby improving (\ref{eq:lb}) to:
\begin{equation}\label{eq:lb2}
g(S) \ge \frac{1}{(d+1)^{d}}\binom{n}{d+1}+O(n^d)
\end{equation}
More generally, if we could show that any core point $p$ of a
$d$-dimensional configuration is contained
in at least $\mu(d)$ simplices, then we can improve the constant
in equation (\ref{eq:lb}) by a factor of $\mu(d)$.

%
%
\section{Colourfully Covering the Core}\label{se:covering}
This leads us to ask: What is the minimum number $\mu(d)$ of simplices that
can contain a core point $p$ in a colourful configuration?  
Given a colourful configuration $\S$ with colourful
sets $S_1, \ldots, S_r$ we can define:
\begin{equation}\label{eq:m}
\m(\S)= \min_{p \in \core(\S)} \csd_\S(p)
\end{equation} 
We remark that if $\core(\S)$ has a non-empty interior, the
minimum in (\ref{eq:m}) will be attained on an open set of
points that are in general position relative to $\S$.

In this notation, our objective is to find the minimum value of $\m(\S)$
over all configurations $\S$ with full-dimensional core in dimension $d$.  
For a fixed $d$, it is clear that some configuration with
$(d+1)$ points in each of $(d+1)$ colours attains this minimum,
which depends only on the dimension.  
Hence we can define:
\begin{equation}\label{eq:mu}
\mu(d)= \min_{\text{$d$ configurations } \S, ~ p \in \core(\S)} \csd_\S(p) 
\end{equation}

One might suppose that $\m(\S)$ is often large.
As a thought experiment, consider choosing a configuration
at random.  If we take $(d+1)$ points in 
${\R}^d$ from a distribution that is ``nice'' and centrally 
symmetric about the origin \zero, the probability that $\zero$ is 
contained in their convex hull is $\frac{1}{2^d}$ (see e.g.~\cite{WW01}). 
This suggests that for random $\S$, a typical value for 
$\csd_\S(\zero)$ would be $\frac{1}{2^d}(d+1)^{d+1}$.  
For a set $S$ of $(d+1)^2$
points in the plane, plugging this value into \bara's analysis
gives us an estimate of $g(S)$ very close to \bara's 
upper bound (\ref{eq:ub}).  However, it is not immediately clear 
if we should expect $\m(\S)$ to be much smaller than $\csd_\S(\zero)$.

If we take a configuration $\S^\triangle$ with 
$S^\triangle_1$ given by $(d+1)$ points
in general position and  
$S^\triangle_1 = S^\triangle_2 = \ldots = S^\triangle_{d+1}$
we get $\m(\S^\triangle) = (d+1)!$.  In Section~\ref{se:cons} we 
exhibit a configuration $\S^\minus$ with $\m(\S^\minus) \le d^2+1$.

In the remainder of the paper, except where noted, we consider 
configurations with $(d+1)$ points in each of $(d+1)$ colours.

\subsection{Preliminaries}\label{se:prelim}
In \cite{BO97}, \bara~and Onn consider the problem of 
{\it colourful linear programming}.  This is the 
algorithmic version of the colourful \cara~ problem:
Given a core point $p$, how can we {\it find} a
colourful simplex containing $p$?
They begin with some preprocessing which is also helpful here.

Take a colourful configuration $\S$ of $(d+1)$ colourful sets
in $\R^d$, $\S = \{S_1, \ldots, S_{d+1}\}$.  
Take $p \in \interior(\core(\S))$.
Without loss of generality we assume that the core point $p=\zero$. 
Given any finite set of points $T \subseteq \R^d$, scaling the points of 
$T$ does not affect whether $\zero$ lies in the convex hull
of $T$ since the coefficients in a convex combination can
themselves be rescaled.
This allows us to normalize 
$\S$ by rescaling its points to unit vectors.

Let $\conv(T)$ be the convex hull of the points in $T$ and 
$\cone(T)$ be the set of non-negative linear combinations 
of points of $T$.  A cone is {\it simplicial} if it 
can be generated by a set of $d$ linearly independent points
in $\R^d$.
If $T \subseteq \R^d$ is a set of points, $\zero \notin T$, but 
$\zero \in \conv(T)$, then $\cone(T)$ must contain a non-trivial
linear subspace of $\R^d$.  A closed, convex cone is called
{\it pointed} if it does not contain such a subspace, so we
summarize this as:
\begin{lemma}\label{le:cones}
Given any finite set of non-zero points $T \subseteq \R^d$, 
$\zero$ is in $\conv(T)$ if and
only if $\cone(T)$ is not pointed.
\end{lemma}
When $T$ is a finite set of points on the unit $d$-sphere 
$\Sph^d \subseteq \R^d$,
Lemma~\ref{le:cones} is equivalent to saying that $\zero \in \conv(T)$ 
if and only if $T$ is not contained in any open
hemisphere of $\Sph^d$.  One direction is proved by building a hemisphere
from a hyperplane through $\zero$ whose normal lies in the interior of
$\cone(T)$ when this cone is pointed.  The other direction is
proved by observing that an open hemisphere never contains both
a point $p$ and its antipode $-p$. 

We would like to put Lemma~\ref{le:cones} in a form that is convenient
for counting how many simplices generated from $T$ contain \zero.  
To do this, we find it helpful to think about the 
antipode of one of the points.
\begin{lemma}\label{le:antipodes}
If $T=\{p_1, p_2, \ldots, p_{d+1}\}$ is a set of non-zero affinely
independent points in $\R^d$,
$\zero$ is in $\conv(T)$ if and only if the antipode 
$-p_{d+1}$ is in $\cone(p_1, p_2, \ldots, p_d)$.
\end{lemma}
\begin{proof}
Let $K = \cone(p_1, p_2, \ldots, p_d)$.
Since $K$ is a cone generated by $d$ linearly independent points in 
$\R^d$, $K$ is simplicial and hence pointed.  
If $-p_{d+1} \in K$, then we can write it as a conic combination
of the remaining $p_i$, that is:
${-p_{d+1} = \sum_{i=1}^d a_i p_i}$
for some $a_1, \ldots, a_d \ge 0$.  Moving the $p_{d+1}$ term 
to the right hand side of the equation and dividing by $1+\sum_{i=1}^d a_i$
gives $\zero$ as a convex combination of the $p_i$'s.
If $-p_{d+1}$ is not in $K$, then we can
strictly separate $-p_{d+1}$ from $K$ with a hyperplane $H$ through \zero.
Then both $K$ and $p_{d+1}$ lie strictly on the same side of $H$, and the
cone generated by $T$ must be pointed.
\end{proof}

\subsection{A Variational Approach}\label{se:var}
Take a point $p$ from a finite set $S \in \R^d$. 
Call a simplex generated by points in $S$ a {\it $p$-simplex} if
$p$ is one of the points used to generate the simplex, and
call a $p$-simplex {\it zero-containing} if it contains $\zero$ in its interior.
Define $z_S(p)$ to be the number of zero-containing $p$-simplices for a
given $S$.

Lemma~\ref{le:antipodes} tells us that $z_S(p)$ is the number of 
simplicial cones generated by $S \setminus \{p\}$ that contain $-p$.  
We find it useful to think about what happens to $z_S(p)$
if we move $p$ while fixing the remaining points of $S$.  
This is particularly illustrative if we confine $p$ to the surface 
of the unit sphere $\Sph^{d}$ centred at \zero.

Let $U = S \setminus \{p\}$ with $|U| = u$.
Initially $z_S(p)$ will be the number of simplicial cones generated
by sets of $d$ points from $U$ that contain $-p$.  Now consider
what happens as $p$ (and hence $-p$) move.  The value of $z_S(p)$ will
stay fixed until $-p$ crosses the boundary of some simplicial cone
from $U$.  These boundaries are defined by the hyperplanes 
generated by $\zero$ and sets of $(d-1)$ points from $U$.  
Taking all $(d-1)$ sets from $U$, we can generate all such boundaries.
They divide the surface of $\Sph^d$ into open {\it cells} that are
$(d-1)$-dimensional open sets.  We can define the {\it depth} of a cell 
of $S$ to be the number of simplicial cones generated by $S$ containing
any given point in the interior of the cell.

Consider moving $p$ along the surface of $\Sph^d$ to a new point
$p'$.  If $-p$ and $-p'$ are in the same cell, we will have 
$z_S(p)=z_S(p')$.  Now suppose $-p$ is in a cell $C$ adjacent to 
the cell $C'$ containing $-p'$.
Then as we move from $-p$ to $-p'$ we cross a single hyperplane $H$
defined by a set $U^0$ of $(d-1)$ points from $U$ belonging to $H$.  
Let's say that 
$-p$ is on the left of $H$ and $-p'$ is on the right.  
For the moment we will assume that only $(d-1)$ points of $U$
lie on $H$.
Let $U^-$ be the
set of $k$ points from $U$ on the left of $H$, and $U^+$ be the
$u-k-(d-1)$ points from $U$ on the right.  Since $-p$ is in
a cell bordered by $H$, it lies in the cone defined by the points
from $U^0$ and $x$ for any point $x \in U^-$.  On the other hand,
$-p$ is separated by $H$ from the cones formed by $U^0$ and $y$ 
for any $y \in U^+$.  Hence $-p$ is contained in exactly $k$
simplicial cones from $S$ generated by $U^0$ and a single other
point.  Similarly, $-p'$ is contained in exactly $u-k-(d-1)$ such
cones.  Simplicial cones that do not contain $U^0$ in their generating 
set will not have $H$ as a facet, so they will contain $-p$ if and only
if they contain $-p'$.  Suppose $-p$ is in $l$ such cones.
Then $z_S(p)=l+k$, while $z_S(p')=l+u-k-(d-1)$.

We conclude that given the value of $z_S(p)$ at some point $p$,
we can in principle compute $z_S(p')$ for any other point $p'$
by tracing a path from $-p$ to $-p'$, and seeing how each hyperplane 
generated from points in $U = S \setminus \{p\}$ divides the points
of $U$.  To do this formally, we need a topological lemma that 
says we can always draw a path between two points on $\Sph^d$ that
crosses only hyperplanes from $U$ (as opposed to passing through
cones generated by fewer than $(d-1)$ points).  
This reduces to the following fact which can be proved using
algebraic topology, see for example \cite{Mun84}:  
\begin{lemma}\label{le:top}
The sphere $\Sph^d$, a $(d-1)$ dimensional manifold, remains path 
connected after
removing finitely many $(d-3)$-dimensional manifolds. 
\end{lemma}

\subsection{Parity}\label{se:parity}
The variational approach to computing $z_S(p)$ explains the
following parity phenomenon:
\begin{proposition}\label{pr:mincfg}
For any colourful configuration $\S$ of $(d+1)$ points in each of $(d+1)$ 
colours in {\rm odd} dimension $d$ and any point $p$ with $\S$ and $p$
in general position, the colourful simplicial depth of $p$ 
with respect to $\S$ is {\rm even}.
\end{proposition}
The authors were surprised by this fact while experimenting with
configurations.  However, it is easy to explain this
via a colourful version of the method described in
Section~\ref{se:var}.  Suppose we begin with a 
configuration $\S^0$ with $(d+1)$ points in each of $(d+1)$
colours clustered near the North Pole of $\Sph^d$.  (We remarked
in Section~\ref{se:prelim} that it is sufficient to consider 
configurations on the surface of $\Sph^d)$.  
Then we can move one point at a time from its initial position
in $\S^0$ to its final position in $\S$ generating a sequence of
configurations $\S^0, \S^1, \S^2, \ldots, \S^{(d+1)^2}=\S$.
Clearly $\csd_{\S^0}(\zero)=0$.  As we move a given point $p_i$ of colour
$j$ from its
initial position in $\S^0$ (and $\S^i$) to its final position
in $\S$ (and $\S^{i+1}$), we need only to know what happens
when the antipode $-p_i$ crosses colourful hyperplanes defined by 
a set of $(d-1)$ points of $(d-1)$ colours, and not of colour $j$.
Such a colourful hyperplane $H$ will miss only one other colour, $j'$.
There will be $k$ points of colour $j'$ on one side of $H$, 
and $(d+1-k)$ on the other side.  
Here we are assuming that the points from $\S$ are in general position, 
but we can argue by continuity that this assumption is not necessary.
As $-p_i$ crosses
$H$ the number of simplicial cones containing $-p_i$ generated by points 
from $H$ and a point of colour $j'$ changes from $k$ to $(d+1-k)$.
As long as $(d+1)$ is even, the parity doesn't change. 

Examining this proof, we can see that Proposition~\ref{pr:mincfg} can
be generalized:
\begin{theorem}\label{th:cparity}
If $\:\S=\{S_1, S_2, \ldots, S_r\}$ is a $d$-dimensional colourful 
configuration of points
and for each $i=1,2, \ldots, r$ we have $|S_i|$ even, and $p$ is
and any point $p$ with $\S$ and $p$ in general position, then the colourful
simplicial depth of $p$ with respect
to $\S$ is also even.
\end{theorem}

For monochrome depth, as we move point $p$ around $\Sph^d$
we need to consider all possible hyperplanes formed from 
$S \setminus \{p\}$.  Using the same reasoning as 
Theorem~\ref{th:cparity} we get:
\begin{theorem}\label{th:mparity}
If $S$ is a set of $n$ points in $\R^d$, and $n-d$ is even,
and $p$ is a point such that $S \cup \{p\}$ is in general position, 
then the simplicial depth of $p$ with respect to $S$ is even.
\end{theorem}

\begin{remark}
The variational approach suggested in Section~\ref{se:var}
has appeared in various guises in discussions of monochrome
simplicial depth.  In particular, it underlies many of the
algorithms suggested for computing monochrome simplicial depth.
Several such algorithms have been proposed recently, see for
example the discussion in \cite{Alo05}.  
Many of these focus on the 2-dimensional problem, 
\cite{GSW92} and \cite{CO01} use variational ideas in
3- and 4-dimensional algorithms.

For this reason, 
we believed that Theorem~\ref{th:mparity} existed as
folklore for some time.  
Baker remarks on the two-dimensional version in a
recreational mathematics note \cite{Bak78}, but this fact,
which impressed the authors with its simple elegance, 
seems curiously neglected in the literature.  
We speculate that one reason for this
is that in statistics the focus has been on computing the depth of the 
sample points themselves, which are not in general position and 
do not retain nice parity conditions.
\end{remark}

\subsection{Configurations with Small Minimal Colourful Depth}\label{se:cons}
We now describe how to build a colourful configuration $\S^\minus$
that contains $\zero$ in its core, but where only $d^2+1$ colourful
simplices contain \zero.  Our strategy is to fix the first $d$
colourful sets $S_1^\minus, S_2^\minus, \ldots, S_d^\minus$ and then consider possible
placements of $(d+1)$ points $p_1, p_2, \ldots, p_{d+1}$ to form
$S_{d+1}^\minus$.  We will place the points from 
$S_1^\minus \cup S_2^\minus \cup \ldots \cup S_d^\minus$ on
the sphere $\Sph^d$ in such a way that some regions of $\Sph^d$
are sparsely covered by simplices from 
$S_1^\minus \cup S_2^\minus \cup \ldots \cup S_d^\minus$.

We begin by fixing $\epsilon = \frac{1}{100d}$.  
We will place the points from $\S^\minus$ in three locations on $\Sph^d$.  
The first on the {\it Tropic of Capricorn}, which we define to be
the set of points on $\Sph^d$ whose $d$th coordinate is $-2\epsilon$.  
The second is on the {\it Tropic of Cancer}, whose $d$th coordinate
is $\epsilon$.  
The two tropics are topologically copies of $\Sph^{d-1}$, but
unlike their namesakes they are not equally spaced from the equator.  
The final region is the {\it polar region}
of points in $\Sph^d$ which are within $\epsilon$ of the North Pole
$p_{\north} = (0,0,\ldots,0,1)$.
\begin{figure}[h!bt]
\begin{center}
\includegraphics[height=5.4cm]{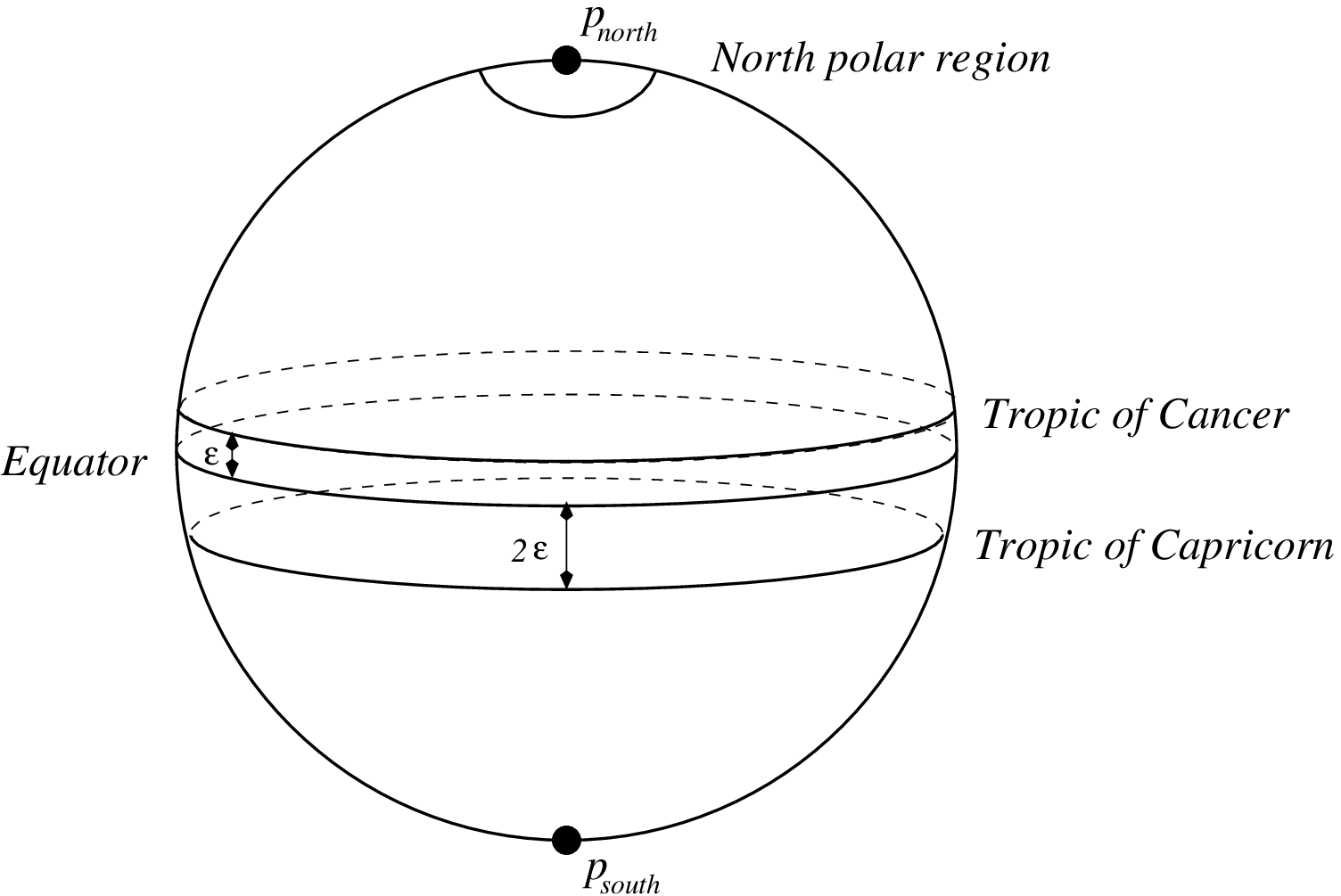} 
\caption{3-dimensional illustration of the regions used in constructing
         $\S^\minus$.}\label{fig:regions}
\end{center}
\end{figure}

Now let's fix the positions of the points 
$\{x_1, x_2, \ldots, x_{d+1}\} \in S_1^\minus$.
Take:
$$x_1 = (\sqrt{1-4\epsilon^2}, 0, 0, \ldots, 0, -2\epsilon)
\qquad \text{ and } \qquad 
x_2 = (-\sqrt{1-\epsilon^2}, 0, 0, \ldots, 0, \epsilon)$$
Note that the line segment between $x_1$ and $x_2$ passes just
below the origin in the sense that it contains a point whose 
first $(d-1)$ coordinates are 0, and whose $d$th coordinate is
negative (and small).
We now place the remaining points $x_3, \ldots, x_{d+1}$ in the
polar region in such a way as to ensure that 
$\zero \in \interior(\conv(S_1^\minus))$.
For $d=2$ we can do this by placing $x_3$ at the North Pole.
For $d \ge 3$ we can place the points on the section of the 
Arctic Circle (points with distance $\epsilon$ to the North Pole)
with zero initial coordinate.  Topologically the Arctic Circle is a copy of
$\Sph^{d-2}$; we can take $x_3, \ldots, x_{d+1}$ to be the vertices
of a regular simplex inscribed on this sphere.

The points of colours $2,3, \ldots, d$ are chosen similarly.
The first points from each of the $d$ colours are arranged in a regular
simplex on Capricorn.  The remaining points in the
same relative position to the first point, so that each $S_i^\minus$
is a rotation of $S_1^\minus$ around the $d$th coordinate axis.
In particular, for each $i=1,2, \ldots, d$, the second point 
of $S_i^\minus$ will lie on Cancer and the final $(d-1)$ points will
lie in the polar region.

We finish our construction by considering possible placements
of the points $p_1, \ldots, p_{d+1}$ of $S_{d+1}^\minus$.
We want to place the $p_i$'s in such a way that their antipodes
(the $-p_i$'s) are contained in few colourful simplicial cones generated
from $\S^\minus$.  

Consider the cell $C_{\south}$ defined by colours $1,\ldots,d$ of
$\S^\minus$ on $\Sph^d$ which contains the South Pole 
$p_{\south}=(0,0, \ldots, 0,-1)$.
We claim this is exactly the intersection of $\Sph^d$ with 
the single colourful simplicial cone $K_{\Cap}$ defined by the $d$ 
colourful points on Capricorn.  
This follows since any other colourful cone 
is generated by a set of $d$ coloured points 
chosen from Capricorn, Cancer and the northern polar region.
Fix such a cone and call these sets ${G_{\Cap}}$, $G_{\Can}$ and
$G_{\Pole}$ and let $K_G = \cone(G_{\Cap} \cup G_{\Can} \cup G_{\Pole})$.  
We assume that we have $|G_{\Cap}| < d$.
We need to show that 
$\interior(K_{\Cap}) \cap \interior(K_G) = \emptyset$.
To do this, we find a hyperplane separating $K_{\Cap}$ and $K_G$.
If $G_{\Cap}=\emptyset$ the hyperplane through the Equator will do.
Otherwise, take the colours from $G_{\Cap}$ and consider any
facet $F$ of $K_{\Cap}$ containing generators of each of these colours.  
Then $F$ separates $K_{\Cap}$ from all the polar points and all the 
Cancer points of colours from $\{1,2,\ldots,n\} \setminus G_{\Cap}$.
(To be absolutely
proper, in higher dimension we would have to move Capricorn up 
towards the equator to ensure the separation of the Cancer points, 
i.e. we would have to reduce the constant $2\epsilon$ to 
$(1+\delta)\epsilon$ for some $\delta>0$.)
This completes the proof.  We conclude that the cell $C_{\south}$ 
is covered only by the colourful cone $K_{\Cap}$ and 
closely approximates the spherical cap bounded by Capricorn.

It is a good strategy to place the antipodes $-p_i$ in $C_{\south}$.
If we do this for all of $S_{d+1}^\minus$, however, 
the resulting configuration will not have $\zero \in \conv(S_{d+1}^\minus)$
($S_{d+1}^\minus$ would certainly be contained in an open hemisphere).
So we must have at least one antipode, say $-p_1$ above Capricorn.
Indeed, if we place the remaining $-p_i$ below Capricorn, we
would need to have $-p_1$ above the ring of the antipodes
of Capricorn.  More precisely, this is the set of points on $\Sph^d$ 
with final coordinate value exactly $2\epsilon$.  In particular,
it is above Cancer.

Let $A=\{a_1, a_2, \ldots, a_d\}$ be the points from 
$S_1^\minus, S_2^\minus, \ldots, S_d^\minus$ 
on Capricorn.  Similarly, let $B=\{b_1, b_2, \ldots, b_d\}$ be 
the points on Cancer.  Let's count how many simplicial cones
from $\S^\minus$ must contain $-p_1$ if we place $-p_1$ above Cancer.
To do this, we start with $-p_1$ in $C_{\south}$ and then move it
above Cancer noting which cell boundaries it crosses as suggested
in Section~\ref{se:var}.  This structure of the cell boundaries
is a topological question, so we find it convenient to remove the
$p_{\south}$ and equate $\Sph^d$ with $\R^{d-1}$.  

With the exception of the single colourful cone that contains
$C_{\south}$, the colourful simplicial cones generated by $\S^\minus$
correspond to colourful simplices in $\R^{d-1}$.  
The polar points on $\Sph^d$ will be clustered near the origin
in $\R^{d-1}$.  Let $A'=\{a_1', \ldots, a_d'\}$ and 
$B'=\{b_1', \ldots, b_d'\}$ be the projections of $A$ and $B$
respectively in $\R^{d-1}$.  
Then $\conv(A')$ and $\conv(B')$ form nested simplices which contain
the projection of the polar region.   The boundaries of the colourful 
simplicial cones on $\Sph^d$ map to facets of simplices in $\R^{d-1}$;
both are defined by sets of $(d-1)$ colourful points.
Moving $-p_1$ from below Capricorn to above
Cancer corresponds to moving $-p_1'$ from outside $\conv(A')$ to
inside $\conv(B')$.

Let us now see what simplicial facets $-p_1'$ must cross to do this.
If we keep $-p_1'$ far away from the $a_i'$ and $b_i'$'s themselves,
we can avoid any facets involving the polar points: These facets
involve at most $(d-2)$ generators from $A'$ and $B'$, and hence have
ends that are at most $(d-3)$ dimensional manifolds in 
$\conv(A') \setminus \interior(\conv(B'))$.  The ends can be avoided by
Lemma~\ref{le:top}. 

This still leaves $d 2^{d-1}$ colourful facets defined by choosing
$(d-1)$ colourful points from $A'$ and $B'$.  We can enumerate them by
first choosing an index (colour) to omit and then representing the
choices of $a_i'$'s and $b_i'$'s by a 0-1 vector of length $(d-1)$.  
Letting 0 represent the choice of an $a_i'$, $\conv(A')$ is bounded by
the facets defined by $d$ index choices and a vector of 0's, while
$\conv(B')$ is bounded by the facets defined by $d$ index choices 
and a vector of 1's.  In fact there are $2^d$ colourful simplices 
generated by $A'$ and $B'$, and they are enumerated by 0-1 vectors of 
length $d$.  Their facets are enumerated by choosing an
index to drop from the enumerating sequence.  Therefore the sums
of the 0-1 vectors enumerating the facets of a given simplex
can differ by at most 1.

Now start with $-p_1'$ outside $\conv(A')$.  
To bring $-p_1'$ inside $\conv(B')$, we must start by bringing it
into $\conv(A')$.  This involves crossing some boundary face of
$\conv(A')$, say the one defined by $a_1, \ldots, a_{d-1}$.  
This is enumerated as $(d,0,0,\ldots,0,0)$.  We can proceed through
facets $(d-1,0,0,\ldots,0,1)$, $(d-2,0,0,\ldots,0,1,1)$
until finally we cross $(1,1,1, \ldots, 1)$ into a cell of $\conv(B')$.
This involves $d$ facet crossings, which is minimal since at
each crossing we can only add a single 1 to the 0-1 part
of the enumerating vector.

We claim that as $-p_1'$ crosses each facet, it makes a net 
gain of $d-1$ containing simplices.  
At the first facet, $(d,0,0,\ldots,0,0)$, $-p_1'$
leaves the single exterior simplex defined by the points $A'$ 
projected from Capricorn and enters the $d$ simplices defined
by $a_1', \ldots, a_{d-1}'$ and 
the $d$ points of colour $d$ other than $a_d'$.  
At subsequent facet crossings, the same thing happens for the remaining
colours: $-p_1'$ leaves the simplex defined by the crossing
facet and $a_i'$.  As $-p_1'$ leaves, it enters the simplices
defined by this facet and the $d$ remaining points of colour $i$.
Hence the number of simplices containing $-p_1'$ immediately after 
crossing into $\conv(B')$ is exactly $1+d(d-1)$.  

We will now return our attention to $\Sph^d$.
Denote by $C_p$ the cell containing $-p_1$ whose projection lies
inside $\conv(B')$.  From our construction, $C_p$ is a cell above
Cancer.  We want to claim that in fact it contains some point above 
the set of antipodes of Capricorn, that is, a point whose antipode
is in $C_{\south}$.
This is a complicated geometric calculation.  However, we 
observe that nothing in our topological argument above changes if we
change the constant $2\epsilon$ in our definition of Capricorn to 
$c \: \epsilon$ for any $c \ge 0$.  In particular, the cell $C_p$
does not degenerate if we move the antipodes of Capricorn
towards Cancer by decreasing $c$ to 1.  
Therefore for some $c>1$ (this condition maintains 
$\zero \in \interior(\conv(S_i^\minus))$ for $i=1,\ldots,d$), 
$C_p$ includes some point 
above the antipodes of Capricorn.  Any such $c$ and
point in $C_p$ would be sufficient for our construction.
We have used $c=2$ for concreteness and take it as an article
of faith that this is a small enough for our choice of $\epsilon$.
\begin{figure}[h!bt]
\begin{center}
\includegraphics[height=5.4cm]{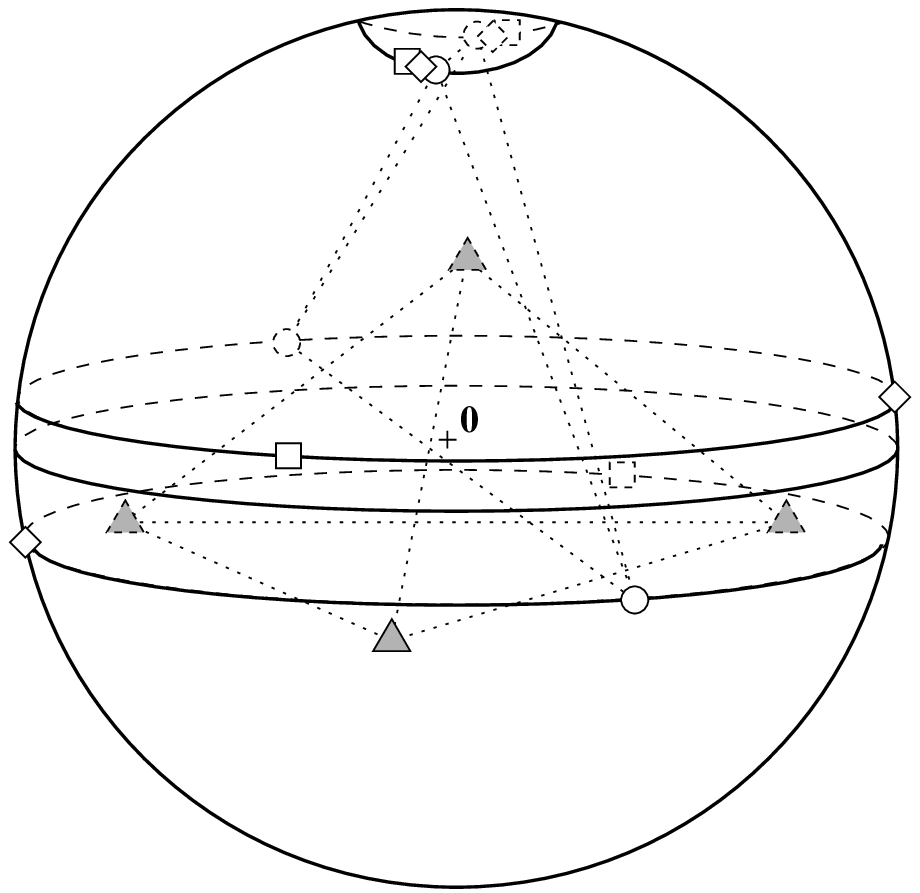} 
\caption{Placement of points of colours
           $1,2,3$ and antipodes of colour $4$ in 
           the 3-dimensional $\S^\minus$.}\label{fig:placement}
\end{center}
\end{figure}

The construction can now be completed.
Take $-p_2$ to be the midpoint of shortest spherical segment between 
Capricorn and $p_1$ (which lies below Capricorn).  
Let $z < -2\epsilon$ be the final coordinate of $-p_2$ and
arrange the remaining points so that $-p_2, -p_3, \ldots, -p_{d+1}$
form a regular simplex on $\Sph^d \cap \{x \in \R^d | x_d=z\}$.
Then $\zero$ is in the convex hull of the $-p_i$ (and hence $S_{d+1}^\minus$).
Finally we can calculate $\csd_{\S^\minus}(\zero)$
from the location of the $-p_i$: $\zero$ lies in $1+d(d-1)$ colourful
simplices generated with $p_1$ and one simplex each including
$p_2, p_3, \ldots, p_{d+1}$.  Hence:
$$\csd_{\S^\minus}(\zero)=1+d(d-1)+d = d^2+1$$

\begin{remark}\label{re:alt}
There are other nice configurations with
$\csd_{\S^\minus}(\zero)=d^2+1$.
Consider a configuration $\S'$ similar to $\S^\minus$ but with the tropics
pushed to the north, taking Cancer's final coordinate to $3\epsilon$
and Capricorn's to $-\epsilon$.  We can then move each of
$-p_1, \ldots, -p_d$ across Capricorn and the equator through a
single boundary facet.  Finally place $-p_{d+1}$ at the South Pole.
Using the same analysis as above, we have $p_1, \ldots, p_d$ points 
forming $1+(d-1)$ simplices containing \zero, and $p_{d+1}$
forming one such simplex for a total of $d^2+1$.

Both $\S^\minus$ and $\S'$ have symmetry for the first $d$ colours,
but not the last one.  We can also propose a configuration $\T$
with symmetry between all the colours.  Follow the recipe for
$\S^\minus$ but place one point of each colour on Cancer and Capricorn 
and place the remaining points in the polar region.  
This brings a number of technical
difficulties, however.  The points will not be in general
position, since the tropical hyperplanes include $(d+1)$ points.
It is also a bit less natural to evenly space $(d+1)$ points
on copies of $\Sph^{d-1}$, indeed for $d=2$ this construction
does not make sense.  When there is a nice way to do this for $d \ge 3$
this (e.g.~4 points on $\Sph^2$) we may end up with some points
being antipodes.  This would cause $\zero$ to be on the faces of
some simplices and increase its colourful simplicial depth.  
Most of these problems can be fixed by perturbing $\T$, but even
so $\T$ is not well-suited to our proof technique.
One might also consider configurations that are not confined
to the sphere.
\end{remark}

\subsection{Evaluating $\mu(d)$}\label{se:min}
The configuration $\S^\minus$ of Section~\ref{se:cons} satisfies
$\m(\S^\minus) \le d^2+1$ where $\m(\S)$ is the minimum colourful simplicial
depth of core point defined in (\ref{eq:m}).
We would like to prove that $\m(\S^\minus)=d^2+1$ and in fact
that for any colourful configuration $\S$ we will have
$\m(\S) \ge d^2+1$, or equivalently $\mu(d) \ge d^2+1$.  
The second half of this proposition clearly
implies the first.  We suggest it is also more approachable since 
we can move the core point of minimum depth to $\zero$ during preprocessing,
whereas a direct attack on $\m(\S^\minus)$ requires understanding the
shape of the core of $\S^\minus$.

\bara's original Colourful \cara~ theorem is exactly
that $\mu(d) \ge 1$.  He further shows that 
for any $\S$ any colourful point from $\S$ is part of some
generating set for a colourful simplex containing \zero.  This
immediately yields $\mu(d) \ge d+1$.  In $\S^\minus$ 
we see that $p_2, p_3, \ldots, p_{d+1}$ all generate a 
{\it unique} colourful simplex containing $\zero$.  Thus
the minimum number of colourful simplices containing 
$\zero$ generated by an arbitrary point in a configuration is 1.
To get a stronger lower bound than $\mu(d) \ge d+1$ we need to 
understand some global information about configurations.

\begin{lemma}\label{le:depthd}
Fix the sets $S_1, \ldots, S_d$ from a colourful configuration $\S$
with $\zero$ in its core, and consider the cells created
on $\Sph^d$ by the colourful simplicial cones from these sets.
Then every cell has depth at least 1, and if there is a cell of
depth 1 it is unique and all other cells have depth at least $d$.
\end{lemma}
\begin{proof}
The fact that every cell has depth at least 1 is equivalent to the 
fact that every colourful point generates some colourful simplex
that contains \zero, proved in \cite{Bar82}.
Suppose now that there is a cell $C$ of depth 1.
Any point exiting $C$ through a bounding hyperplane will be
exiting some colourful simplex.  Since the depth of $C$ is 1,
this will always be the same simplex.
Thus the extreme points of $C$ must be a colourful set
$A=\{a_1, \ldots, a_d\}$ with $a_i \in S_i$ generating this simplex.  
We can puncture $\Sph^d$ at $p \in C$ and project $\Sph^d \setminus \{p\}$ 
into $\R^{d-1}$.  The $a_i$'s project to a set 
$A'=\{a_1', \ldots, a_d'\}$ that forms a $(d-1)$-simplex 
in $\R^{d-1}$.  
The remaining colourful points project to points in $\conv(A')$.

Take a point $q$ inside $\conv(A')$.  
We want to show that $q$ is contained in at least $d$ 
colourful simplices in addition to $\conv(A')$ after projection.
To do this, it is sufficient to show that if we take any colourful
set $B'=\{b_1', \ldots, b_d'\}$ of projected points with $b'_i$ of
colour $i$ and $A' \cap B' = \emptyset$, then $q$ is in some colourful
simplex generated from points of $A' \cup B'$ with some generators
from $B'$.  Equivalently, we want to show that $\conv(A')$ is
covered by colourful simplices generated from $A' \cup B'$
(excluding $\conv(A')$ itself from the covering).
Then by partitioning the projections of the colourful
points into $(d+1)$ colourful sets $A', B_1', B_2', \ldots, B_d'$
we get Lemma~\ref{le:depthd}.

Consider the collection $\X$ of colourful $(d-1)$-simplices 
generated by $A'$ and $B'$ in $\R^{d-1}$ and let 
$\tilde{X}$ be the set of points 
contained in the colourful simplices of $X$ other than $\conv(A')$.
The elements of $\X$ are all the simplices
formed by taking for each colour $i=1,2,\ldots,d$ either 
$a_i'$ or $b_i'$ as a generating vertex.  
This construction resembles the $d$-dimensional 
cross-polytope $\beta_d$ (the dual of the $d$-cube), a regular polytope
in $\R^d$ with $2d$ vertices and $2^d$ facets.  
The cross-polytope $\beta_d$ is generated by taking as vertices the
standard unit vectors $E^+ = \{e_1, \ldots, e_d\}$ 
and their negatives $E^- = \{-e_1, \ldots, -e_d\}$.
The facets of $\beta_d$ are the convex hulls
generated by choosing for each $i=1,\ldots,d$ either $e_i$ or $-e_i$.

We can see that $\X$ is obtained from $\beta_d$ as follows:
We have $A' \cup B' \subset \R^{d-1}$.
Embed $\R^{d-1}$ as an affine subspace $\aff(A')$ in $\R^d$.  
Take $\P$ to be an affine hyperplane in $\R^d$ parallel to $\aff(A')$.
For $i=1, \ldots, d$ let $p_i$ be the intersection
point of $\P$ with the line through $b_i'$ perpendicular to $\R^{d-1}$.
Let $P =\{ p_1, \ldots p_d\}$ and generate a set $Q$ of 
$(d-1)$-simplices
by taking for each $i = 1, \ldots, d$ either $a_i'$ or $p_i$.
By construction $X$ is the projection of $Q$ into $\aff(A')$.  
Now we claim that $Q$ is a continuous image of the facets of $\beta_d$.
We can exhibit such a map by first finding an affine transformation $T_1$
with $T_1(e_i)=a_i'$ for $i=1,\ldots,d$ and $T_1(\aff(E^-))=\P$.  
Note that $T_1(\beta_d)$ is a polytope.
Then applying a further affine transformation $t_2$ to $\P$ with
$t_2(-e_i)=p_i$ for $i=1, \ldots, d$ and extending this to 
$T_2$ on $\R^d$ so that $T_2$ fixes $\aff(A')$, 
we see that the composition $T_2 \circ T_1$ is the required map.

We proceed by contradiction.  Assume that $\tilde{X}$ 
does not cover $\conv(A')$.  Then we can find a retraction of
$\tilde{X}$ to its boundary $\partial(\conv(A'))$.
By composing $T_2$, the projection taking $Q$ onto $X$ and the
retraction of $\tilde{X}$, we get 
a retraction of $T_1(\beta_d) \setminus \conv(A')$
onto $\partial(\conv(A'))$.  
However, $T_1(\partial(\beta_d))$ is a $d$-dimensional 
polytope topologically equivalent to $\Sph^d$ and hence 
$T_1(\partial(\beta_d)) \setminus \conv(A')$ is topologically equivalent to a
$(d-1)$-dimensional disk $\B^{d-1}$.  But $\partial(\conv(A'))$ is
topologically equivalent to $\Sph^{d-1}$ and a well-known theorem
of algebraic topology says that there does not exist a retraction
of $\B^{d-1}$ to $\Sph^{d-1}$ (see for example section 21 of \cite{Mun84}).  
This is the required contradiction, hence the colourful simplices 
of $\X \setminus \conv(A')$ cover $\conv(A')$.  
\end{proof}

\begin{corollary}\label{co:muge2d}
The minimum colourful simplicial depth of any core point in any
colourful configuration is at least $2d$.
That is, we have $\mu(d) \ge 2d$.
\end{corollary}
\begin{proof}
It suffices to prove this for a configuration $\S$ with $(d+1)$
in $(d+1)$ colours.
Observe that
if we have no cell of depth 1 then each of the $(d+1)$ points
of $S_{d+1}$ will generate at least two colourful simplices
containing \zero, and if we do have such a cell $C$, we must place
at least one point, say $p_1 \in S_{d+1}$ outside of $C$ to 
get $\zero \in \conv(S_{d+1})$.  Then $p_1$ is generates at least $d$ simplices
containing $\zero$ in addition to the $d$ required of the remaining
points in $S_{d+1}$.
\end{proof}

\subsection{The Two-dimensional Case}\label{se:2d}
We will briefly illustrate our methods by describing how core
points can be contained in configurations in $\R^2$.
Consider such a configuration $\S=\{X,Y,Z\}$ with core point $p$.
We assume general position, and as discussed in Section~\ref{se:prelim},
we may without loss of generality take the core point
$p=\zero$ and place the points of $\S$ on the unit circle $\Sph^2$.  

Then the points of $X$ and $Y$ divide $\Sph^2$ into six segments.
Let $X=\{x_1,x_2,x_3\}$, $Y=\{y_1,y_2,y_3\}$.
These points generate 9 simplicial cones and divide $\Sph^2$ 
into 6 segments.  The boundaries between cones are simply the
rays through the $x_i$'s and $y_i$'s.
Because no three points of $X$ or $Y$ lie in the same half-circle,
each hyperplane through $\zero$ and $x_i$ divides
the $y_i$'s 2 to 1 and vice-versa.  
Then as the antipode of a point from $Z$ crosses $x_i$ or $y_i$ 
the number of containing simplicial cones changes by exactly one.

To get a configuration $\S^\minus$ where only 5 simplices contain 0,
we take $x_1=(-\sqrt{1-4\epsilon^2}, -2\epsilon)$, 
$x_2=(\sqrt{1-\epsilon^2},\epsilon)$, 
$x_3=(-\epsilon,\sqrt{1-\epsilon^2})$, 
$y_1=(\sqrt{1-4\epsilon^2},-2\epsilon)$, 
$y_2=(-\sqrt{1-\epsilon^2},\epsilon)$, 
and 
$y_3=(\epsilon,\sqrt{1-\epsilon^2})$.
Observe there is a large cell of depth 1 between $x_1$ and $y_1$.
The reader can verify that the sequence of colourful cell
depths is: 1,2,3,4,3,2. 
\begin{figure}[h!bt]
\begin{center}
\includegraphics[height=1.2cm]{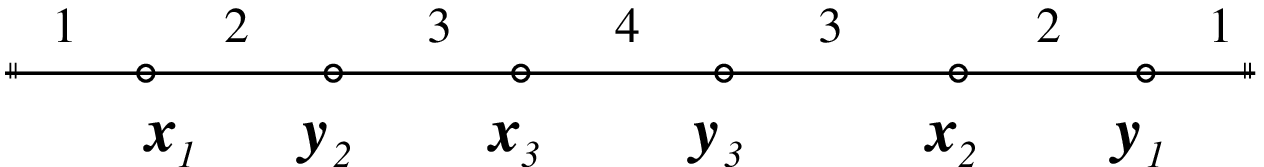} 
\caption{Covering depths for a circle with a depth 1 cell.}\label{fig:circle}
\end{center}
\end{figure}

Let $Z=\{z_1, z_2, z_3\}$. 
Place $z_2=(-\sqrt{1-9\epsilon^2},3\epsilon)$ and 
$z_3=(\sqrt{1-9\epsilon^2},3\epsilon)$ so that their
antipodes lie between $x_1$ and $y_1$.  They each generate one
simplex containing \zero.  Finally, to ensure that $\zero \in \conv(Z)$,
we see that $-z_1$ must lie above $y_2$ and $x_2$.  Take
$z_1=(-\sqrt{1-16\epsilon^2}, -4\epsilon)$.  
Then $-z_1$ is contained in 3 colourful
simplicial cones generated by $X$ and $Y$.  This configuration
has $\zero$ in the interior of its core and $\zero$ lies in $1+1+3=5$ 
colourful simplices.
\begin{figure}[h!bt]
\begin{center}
\includegraphics[height=5cm]{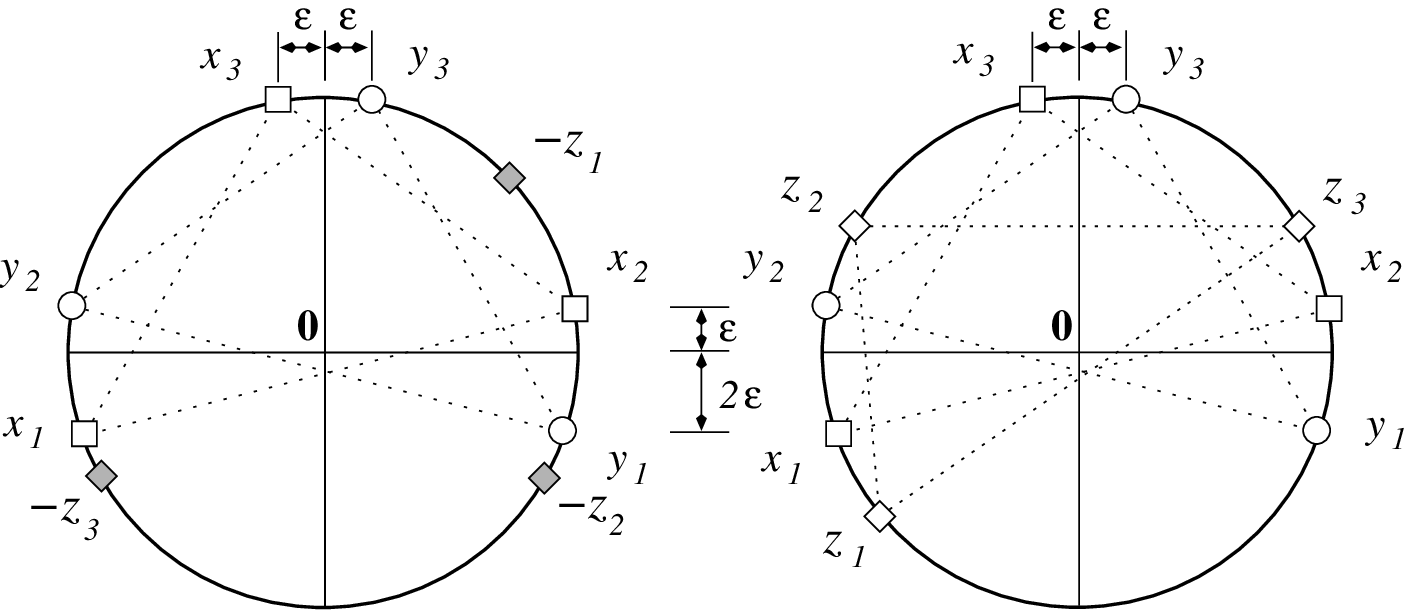} 
\caption{Valid configuration $\S^\minus$ in dimension 2
with $\csd_{\S^\minus}(\zero)=5$.}\label{fig:min2d}
\end{center}
\end{figure}

Using the analysis in Section~\ref{se:min} we see that the
cells generated by $X$ and $Y$ have colourful covering depth at
least 1.  If no cell attains this, then our configuration must
yield at least 6 colourful simplices containing \zero.
If some cell has depth 1, we can place at most two of the
$-z_i$'s in this cell.  The remaining $z_i$ must then have
depth at least 2, for a minimum of 4.  In fact, we can strengthen
this to show that our configuration is minimal  
by observing that we cannot place all of $Z$ in two adjacent cells.
We conclude that $\mu(2)=5$.
A similar observation in three dimensions shows that 
$\mu(3) \ge 8$.  Given the construction of Section~\ref{se:cons}
and Proposition~\ref{pr:mincfg} we know that $\mu(3)$ is
either 8 or 10.
\bara~and Matou{\v s}ek \cite{BM06} have shown that
$\mu(3) \ge 9$ which combined with Proposition~\ref{pr:mincfg} 
implies that
$\mu(3)=10$.

%
%
\section{Conclusions}\label{se:dis}
Let us return to our original goals.  Using the bound $\mu(d) \ge 2d$
from Section~\ref{se:min}, we see that we can improve \bara's
lower bound (\ref{eq:lb2}) for
the depth of the monochrome simplicial median to:
\begin{equation}\label{eq:lb3}
g(S) \ge \frac{2d}{(d+1)^{d+1}}\binom{n}{d+1} + O(n^d)
\end{equation}
This is a modest improvement.  Unfortunately, the construction
in Section~\ref{se:cons} shows that simply bounding $\mu(d)$ 
cannot give a stronger bound than:
\begin{equation}\label{eq:lbc}
g(S) \ge \frac{d^2+1}{(d+1)^{d+1}}\binom{n}{d+1} + O(n^d)
\end{equation}

Quite recently, Wagner proved exactly the bound (\ref{eq:lbc}) in
his thesis \cite{Wag03} as a special case of his
First Selection Lemma.  
This is, to our knowledge, the first improvement of
(\ref{eq:lb2}) since \bara's original paper \cite{Bar82}.
Wagner's result uses a continuous version of the 
Upper Bound Theorem for polytopes and other techniques from
probability without any reference to colouring.  
We find the appearance of the constant $d^2+1$, which for
us arrives from colourful combinatorics, quite remarkable.

\subsection{Bounds for Core Point Depth}
Recalling that $\m(\S)$ is the minimum value of a
core point in a configuration $\S$ and that $\mu(d)$ is the
minimum value of $\m(\S)$ over all $d$-dimensional 
colourful configurations $\S$, our main result is:

\begin{theorem}\label{th:bounds}
The minimal colourful simplicial depth of any interior core point in
any colourful configuration is between $2d$ and $d^2+1$.
That is, we have: $2d \le \mu(d) \le d^2+1$.
\end{theorem}

\begin{conjecture}\label{co:mu}
The minimum colourful simplicial depth of any interior core point in any 
colourful configuration is $d^2+1$.  That is, 
we have $\mu(d) = d^2+1$.
\end{conjecture}
This conjecture implies that the configuration $\S^\minus$ 
minimizes $\m(\S)$ for $d$-dimensional colourful 
configurations.  It would also give an elementary proof of (\ref{eq:lbc}).  
It is easy to see that this holds for $d=1$.  As we noted
in Section~\ref{se:2d}, Conjecture~\ref{co:mu} holds for $d=2$ and $d=3$.
The non-uniqueness of configurations attaining $\m(\S)=d^2+1$
suggests that any such proof cannot be completely trivial
but it may be possible to do this through improved bookkeeping.
The authors generated random low-dimensional configurations 
by computer and did not find any counterexamples to 
Conjecture~\ref{co:mu}.

\begin{remark}\label{re:improved}
The lower bound for $\mu(d)$ was improved very recently 
independently by 
\bara~and Matou{\v s}ek \cite{BM06} and 
Stephen and Thomas \cite{ST06} 
to 
$\displaystyle \max\left(3d, \frac{1}{5}d(d+1)\right)$ for $d>2$
and $\displaystyle \left\lfloor \frac{(d+2)^2}{4} \right\rfloor$ 
respectively.  We know that $\mu(1)=2$, $\mu(2)=5$ and $\mu(3)=10$.
Combining the improved bounds with the parity conditions of
Proposition~\ref{pr:mincfg}
we have the following bounds on $\mu(d)$ for $d > 3$: 
$$12 \le \mu(4)\le 17, ~
  16 \le \mu(5) \le 26, ~ 18 \le \mu(6) \le 37, ~ 22 \le \mu(7) \le 50,$$
and for $d > 7$:
$$
\displaystyle \left\lfloor \frac{(d+2)^2}{4} \right\rfloor 
   \le \mu(d) \le d^2+1.$$
\end{remark}
\vspace{1mm}

It is also natural to ask what type of colourful configuration 
has a core point of {\it maximum} colourful simplicial depth.  
For this question to be interesting, we must fix the number and
size of the colourful sets.  Hence we restrict our attention to
$d$-configurations with $(d+1)$ points in each of $(d+1)$ colours.  
We also require $p$
to lie in the interior of the core since moving to the boundary
of a simplex increases the depth.
We define:
\begin{equation}\label{eq:nu}
\nu(d)= \max_{d \text{ configurations }\S, 
  ~ p \in \interior(\core(\S))} \csd_\S(p) 
\end{equation}

Our method is well suited to analyzing $\nu(d)$ simply
by changing our objective to creating deep cells and placing
antipodes in them.  
We remark that $\nu(1)=2$.  
An analysis similar to that of Section~\ref{se:2d} shows that
$\nu(2)=9$.  The key observation is after placing two sets of
three colourful points on the circle, the sequence of cell depths
that we obtain is either 1,2,3,4,3,2 
or 3,2,3,2,3,2.  In the first case we also need to argue that
the cells of depth at least 3 cover less than half the
circle and that opposite every point of depth 4 is a point of
depth 1.

\begin{figure}[h!tb]
\begin{center}
\includegraphics[height=5cm]{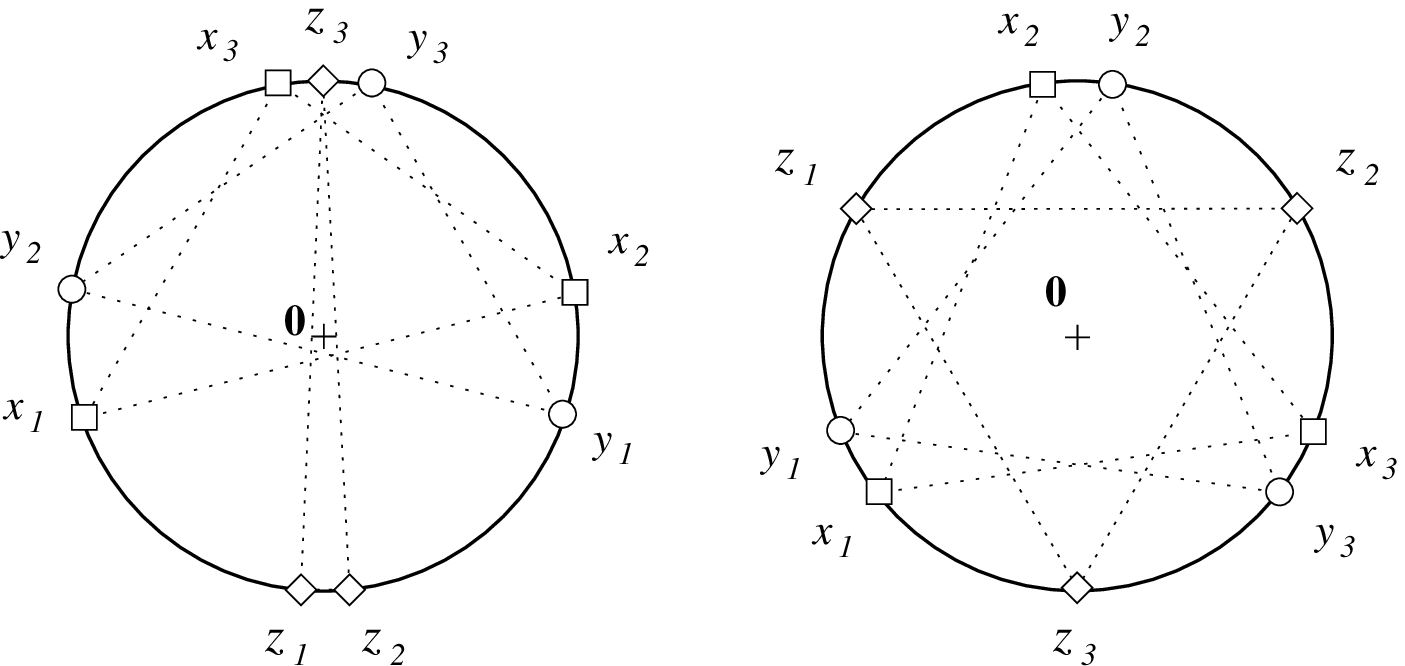} 
\caption{The two configurations in dimension 2
with $\csd_{\S}(\zero)=9$.}\label{fig:max2d}
\end{center}
\end{figure}

The minimal core depth configuration $\S^\minus$
used to prove $\mu(2)=5$ is topologically unique, so it is
interesting to observe that, up to topology, there are two distinct
configurations that contain $\zero$ in 9
colourful simplices.  The first corresponds to the sequence of
cell depths 1,2,3,4,3,2 and contains a point $z_3$
that generates a  unique \zero-containing colourful simplex.
The second corresponds to the sequence 3,2,3,2,3,2 and is a 
combinatorially symmetric configuration where each colourful 
point is in exactly three \zero-containing colourful simplices.
The configurations are illustrated in Figure~\ref{fig:max2d}.

We can build a configuration $\S^+$ with 
$\csd_{\S^+}(\zero) =d^{d+1}+1$ by following the strategy for
$\S^\minus$ but building a deep cell rather than a shallow one.
To do this, we place the polar region points of colour $i$ close
to the geodesic between $p_{\north}$ and the point of colour $i$
on Cancer.  Then $p_{\north}$ is contained in every colourful
cone generated by points from Cancer and the polar region
(in fact these are all the colourful cones containing $p_{\north}$).
Hence the cell $C_{\north}$ containing $p_{\north}$ has depth $d^d$.  
By placing the points of $S^+_{d+1}$ so that $d$ of their antipodes
are in $C_{\north}$ and the final antipode is at $p_{\south}$, 
we get $\S^+$ with $\csd_{\S^+}(\zero) = d \cdot d^d+1$.  
The two-dimensional $\S^+$ appears as the left element of
Figure~\ref{fig:max2d}.

A more symmetric (but similar) construction places one point of
each colour at the vertices of a regular simplex, and the remaining
points surround the antipode of the same colour.

It follows that $\nu(d) \ge d^{d+1}+1$. 
We conjecture that this bound is tight.
As with Conjecture~\ref{co:mu} a computer search did not turn up
any counterexamples.
\begin{conjecture}\label{co:nu}
The maximum colourful simplicial depth of any point in the
interior of the core of any colourful configuration of
$(d+1)$ points in each of $(d+1)$ colours is $d^{d+1}+1$.
That is, we have $\nu(d)=d^{d+1}+1$. 
\end{conjecture}

\begin{remark}\label{re:32}
For any $d$, there exists a colourful configuration $\S$ 
which contains $\zero$ in at least $32\%$ of its colourful simplices.
\end{remark}
A configuration of $(d+1)$ points in each of $(d+1)$ colours 
generates $(d+1)^{d+1}$ colourful
simplices, so Remark~\ref{re:32} follows immediately from 
the construction of $\S^+$.  The minimum fraction of colourful
simplices containing $\zero$ from an $\S^+$ configuration is $82/256$ 
attained when $d=3$.

%
%
\section{Open Questions}
We would like to conclude by mentioning that there are many more
natural questions relating to colourful and monochrome simplicial depth.
The first is:
\begin{question}\label{q:mavg}
What is a typical value of $\m(\S)$ for a random 
configuration $\S$ of $(d+1)$ points in each of $(d+1)$ colours?  
\end{question}
In Section~\ref{se:covering}, we remarked that such random
configurations could be expected to have a simplicial depth
on the order of $\frac{1}{2^d}(d+1)^{d+1}$ at the origin.  
We also gave a colourful configuration 
$\S^\triangle$ that has $\m(\S^\triangle)=(d+1)!$.  
However, $\S^\triangle$ is not in general position.  
Our construction $\S^\minus$ from Section~\ref{se:cons} is 
in general position and has a low value of $\m(\S^\minus)$.
It is not clear if this behaviour is typical, 
i.e.~if most configurations have some point $p$ near the edge of
the core that drags down $\m(\S)$, or if our configuration is
statistically unlikely.  Indeed we can consider the
possibility that {\it all} configurations in general position 
have such a point near the edge of the core.
\begin{question}\label{q:mmax}
What is the maximum value of $\m(\S)$ for a colourful configuration $\S$
of $(d+1)$ points in each of $(d+1)$ colours?  What if $\S$ is not
assumed to be in general position?
\end{question}

We observe that in fact our construction of a colourful
configuration $\S^\minus$ with $\m(\S^\minus)=d^2+1$ contains points of
high colourful simplicial depth, but away from \zero.
This leads us to consider the colourful analogues of the
functions $f(S)$ and $g(S)$ of Section~\ref{se:depth}.
For a colourful configuration $\S$, define:
\begin{equation}\label{eq:fg}
\f(\S)= \max_{p \in \R^d} \csd_\S(p) \qquad \text{ and } \qquad
\g(\S)= \max_{p \text{ in general position}} \csd_\S(p)
\end{equation}
We focus on the case where we have $(d+1)$ colours.
It is clear that given the sizes of the colourful sets 
$S_1,\ldots,S_{d+1}$ comprising $\S$ that the maximum of 
$\f(\S)$ and $\g(\S)$ is $|S_1| \cdot \ldots \cdot |S_{d+1}|$
and is attained by placing the points of each colour at
(or near) the vertices of a simplex.  
If we restrict $\S$ to be a configuration of $(d+1)$ points in 
each of $(d+1)$ colours and take the maximum over the interior of the core, we
get exactly the question of finding $\nu(d)$
(Conjecture~\ref{co:nu}).
We are also interested in lower bounds for $\f(\S)$ and $\g(\S)$.
\begin{question}\label{q:minmax}
For $d$-dimensional configurations consisting of $n$ points in each of 
$(d+1)$ colours, find lower bounds for $\f(\S)$ and $\g(\S)$.
\end{question}

In a survey paper on the colourful \cara~ theorem,
\bara~and Onn \cite{BO97b}
mention that the results of \cite{ABFK92} can be applied to give
a lower bound for $\g(\S)$ when $n$ is large of the form:
\begin{equation}\label{eq:colourlb}
\g(\S) \ge c_d \binom{n}{d+1}
\end{equation}
Unfortunately, the constant $c_d$ is doubly exponential in $d$
so the bound is only non-trivial if $n \gg e^{4d^2}$.  
In particular, it sheds no light on the $n=d+1$ case.

One can also get a lower bound for $\g(\S)$ directly from the
Colourful Tverberg Theorem \cite{ZV92}, which is used to derive the 
results in \cite{ABFK92}:
\begin{equation}\label{eq:tverlb}
\g(\S) \ge \frac{1}{4}\left(\frac{n}{d+1} +3\right)
\end{equation}
This still does not help for $n=d+1$, but 
for small $n$ the bound is stronger than (\ref{eq:colourlb})
and comes with the additional guarantee that colourful
simplices involved are disjoint!  This suggests that there is
much room for improvement.  

\subsection{Monochrome Questions}
The authors would also like to mention that they do not know
the answers to some fairly basic questions about monochrome 
simplicial depth.
Recall the maximum closed and open depth functions $f(S)$ and $g(S)$ 
for a set of points $S$ in $\R^d$  defined in Section~\ref{se:depth}.
\begin{question}\label{q:limits}
Are the points $p$ attaining the maximum $f(S)$ in (\ref{eq:maxf})
always limit points of the set of maxima attaining $g(S)$
in (\ref{eq:maxg})?
\end{question}
We feel that a positive answer to this question would provide
a further natural justification for studying $g(S)$ in place of $f(S)$
when the former is more tractable.  Similarly, it would be interesting
to get conditions on $S$ such that $f(S)$ is not much larger
than $g(S)$.

We are also curious about the expected values of $f(S)$ and $g(S)$:
\begin{question}\label{q:mono}
Given $n$ points in $\R^d$ distributed independently and symmetrically
about \zero, what is the expected deepest simplicial depth of the
resulting configuration?  That is, what is the expected
depth of the simplicial median of the points?
\end{question}
Wagner and Welzl \cite{WW01} give an expression for the 
expected depth of \zero, but $\zero$ will not always
be the deepest point. 
Indeed if $n=d+1$ the expected simplicial depth of $\zero$ will be 
$\frac{1}{2^d}$ while the simplicial median always has depth 1.
For fixed $d$ the expected depth of $\zero$ is
$\frac{1}{2^d}\binom{n}{d+1}$ which has the same asymptotic behaviour as
\bara's sharp upper bound (\ref{eq:ub}) for $g(S)$.
However, when $n$ is not much larger than $(d+1)$, the gap between 
the expected depth of $\zero$ and \bara's upper bound is substantial and it is
not clear to us where the expected depth of the simplicial median lies.

\bara's method of proving 
(\ref{eq:lb}) combined with a solution to Question~\ref{q:mavg}
might lend some insight into Question~\ref{q:mono}, but a direct
approach would be better.

%
%
\section{Acknowledgments}
We would like to thank {\sc Imre B\'ar\'any} for discussions
which triggered this work and the anonymous referees for
suggestions which improved the presentation of this paper.

\bibliographystyle{hacked}

\bibliography{refs}

\end{document}